\documentclass[french,12pt]{amsart}
\usepackage{babel,graphicx}
\newtheorem{theoreme}{Th\'eor\`eme}[section]
\newtheorem{corollaire}[theoreme]{Corollaire}
\newtheorem{lemme}[theoreme]{Lemme}
\newtheorem{proposition}[theoreme]{Proposition}
\theoremstyle{definition}
\newtheorem{definition}[theoreme]{D\'efinition}

\newtheorem{conjecture}[theoreme]{Conjecture}
\theoremstyle{remark}
\newtheorem{remarque}[theoreme]{Remarque}

\newtheorem{question}[theoreme]{Question}
\numberwithin{equation}{section}

\begin{document}

\title[Surfaces branch\'ees et sol\'eno\"{\i}des]{Surfaces branch\'ees et sol\'eno\"{\i}des $\varepsilon$-holomorphes}
\author{Bertrand Deroin}%
\address{Max Planck Institut f\"ur Mathematik -Leipzig}
\email{Bertrand.Deroin@mis.mpg.de}

\thanks{The author acknowledges support from the {\it Swiss National Science Foundation}. }

\begin{abstract} Nous d\'emontrons que pour tout $\varepsilon>0$, 
il existe une lamination compacte par surfaces 
$\varepsilon$-holomorphes dans 
le plan projectif complexe, minimale, et poss\'edant un lacet d'holonomie hyperbolique. 
Nous appelons $\varepsilon$-holomorphe 
une surface r\'eelle $\Sigma$ dont l'angle entre le plan tangent $T\Sigma$ et 
$iT\Sigma$ est uniform\'ement major\'e par $\varepsilon$. 
Ces surfaces sont en particulier symplectiques lorsque $\varepsilon$ est assez petit.  \\ 

\textit{Title.} $\varepsilon$-holomorphic branched surfaces and solenoids. 

\vspace{0.2cm}

\textit{Abstract.} 
We show that for every $\varepsilon>0$, there exists a compact lamination by
$\varepsilon$-holomorphic surfaces in the complex projective plane, minimal, 
and that carries hyperbolic holonomy. We call $\varepsilon$-holomorphic a
real $2$-dimensional surface $\Sigma$ in ${\bf C}P^2$ such that the angle
between $T\Sigma$ and $iT\Sigma$ is uniformly bounded by $\varepsilon$.
When $\varepsilon$ is sufficiently small, such surfaces are in particular symplectic. 

\end{abstract}

\maketitle

\section{Introduction}
Une question encore ouverte est celle de l'existence d'une lamination compacte par courbes holomorphes dans le plan projectif complexe,
ne contenant pas de courbe alg\'ebrique. Notre int\'er\^et pour cette question est la \textit{conjecture du minimal exceptionnel}, qui  
affirme qu'une feuille d'un feuilletage holomorphe singulier du 
plan projectif complexe s'accumule sur une singularit\'e. 
Si une feuille d'un tel feuilletage ne s'accumulait pas sur une singularit\'e,
alors son adh\'erence serait une lamination compacte par courbes holomorphes~; cette lamination serait non triviale dans le sens 
qu'elle ne supporterait pas de mesure transverse invariante, que toutes ses 
feuilles seraient hyperboliques~\cite{Camacho-LinsNeto-Sad1}, et qu'elle 
contiendrait un lacet d'holonomie hyperbolique \cite{Bonatti-Langevin-Moussu}.

\vspace{0.2cm}

Des progr\`es r\'ecents ont \'et\'e faits lorsque l'espace total de la lamination 
est une \textit{hypersurface r\'eelle}. Une hypersurface r\'eelle de ${\bf C}P^2$ est dite 
\textit{Levi-plate} si elle est feuillet\'ee par des courbes holomorphes.
Siu a d\'emontr\'e qu'il n'y a pas d'hypersurface 
Levi-plate compacte de classe $C^8$ 
dans ${\bf C}P^2$~\cite{Siu1}. La r\'egularit\'e 
a \'et\'e am\'elior\'ee par Iordan et 
Cao/Shaw/Wang \`a une hypersurface Levi-plate de classe 
$C^4$ et $C^2$~\cite{Iordan1,C-S-W1}. 
Sous des hypoth\`eses de r\'egularit\'e 
beaucoup plus fortes, et de nature globale, nous avons d\'emontr\'e 
qu'il n'existe pas d'hypersurface Levi-plate compacte \textit{immerg\'ee} dans ${\bf C}P^2$~\cite{Deroin2}.

\vspace{0.2cm}

Le cas d'un \textit{sol\'eno\"{\i}de}, c'est-\`a-dire d'une lamination 
dont l'espace transverse est totalement discontinu, semble bien plus d\'elicat. 
Rappelons qu'une surface lisse de ${\bf C}P^2$ est une courbe holomorphe 
si et seulement si son plan tangent est invariant par multiplication par 
$i=\sqrt{-1}$. \'Etant donn\'e un r\'eel $\varepsilon >0$, appelons 
\textit{$\varepsilon$-holomorphe} une surface lisse ou une famille de surfaces lisses 
du plan projectif complexe dont le plan tangent fait un angle 
inf\'erieur \`a $\varepsilon$ avec son image par multiplication par $i$.
L'angle est mesur\'e avec une m\'etrique hermitienne de 
${\bf C}P^2$, par exemple la m\'etrique de Fubini-Study. 
Dans ce travail, nous montrons que \textit{pour tout $\varepsilon >0$, il existe un 
sol\'eno\"{\i}de compact par surfaces $\varepsilon$-holomorphes 
\textit{plong\'e} de mani\`ere lisse dans le plan projectif complexe, minimal, et 
qui admet un lacet d'holonomie hyperbolique} 
(voir Corollaire \ref{plgt solenoide}). 
En particulier, les feuilles de ces sol\'eno\"{\i}des sont 
symplectiques lorsque $\varepsilon$ est suffisament petit.
De plus, nous conjecturons qu'il existe dans ${\bf C}P^2$ un sol\'eno\"{\i}de 
par courbes holomorphes, minimal, et poss\'edant un lacet d'holonomie hyperbolique (voir Remarque \ref{conjecture}). 

\section{Structure de la d\'emonstration}

Les sol\'eno\"{\i}des par surfaces $\varepsilon$-holomorphes 
que nous obtenons sont inspir\'es d'un sol\'eno\"{\i}de ``abstrait" construit par 
Sullivan pour \'etudier les endomorphismes dilatants du cercle~\cite{Sullivan1}. La
partie \ref{solenoide de Sullivan} est consacr\'ee \`a la construction de ce dernier, et \`a certaines de ses variantes~: 
pour tout $g\geq 0$, nous construisons un sol\'eno\"{\i}de compact $\mathcal S_g$ ``abstrait", 
que nous appelons le \textit{sol\'eno\"{\i}de de Sullivan de genre $g$}. 

\vspace{0.2cm}

Il est possible d'approcher les sol\'eno\"{\i}des $\mathcal S_g$ par des 
\textit{surfaces branch\'ees lisses} $\overline{\Sigma}$. Elles sont obtenues
en partant d'une surface compacte $\Sigma$ avec deux composantes 
de bords $\partial_-\Sigma$ et $\partial _+\Sigma$, que l'on identifie 
par un rev\^etement double
\[ r:\partial_+\Sigma \rightarrow \partial _-\Sigma.\]
Leur structure lisse est explicit\'ee au paragraphe \ref{structure lisse}. 
L'analogue $1$-dimensionel de ce proc\'ed\'e d'approximation est la th\'eorie des \textit{r\'eseaux ferroviaires} 
invent\'ee par Thurston pour \'etudier les laminations g\'eod\'esi\-ques d'une surface hyperbolique~\cite{Bonahon}. 
En dimension sup\'erieure, l'approximation d'un sol\'eno\"{\i}de par des vari\'et\'es branch\'ees a \'et\'e 
\'etudi\'ee par Williams~\cite{Williams} et Gambaudo~\cite{Gambaudo}. 

\vspace{0.2cm}

Au paragraphe \ref{voisinage}, nous \'etudions l'image $\pi (\overline{\Sigma})$ 
d'un plongement lisse $\pi :\overline{\Sigma}\rightarrow V$ dans 
une vari\'et\'e orient\'ee de dimension $4$, et nous donnons des formes locales de la topologie de leur voisinage (Th\'eor\`eme 
\ref{voisinage global}). 
Essentiellement, ils ne d\'ependent que d'un invariant $N(\pi(\overline{\Sigma}),V)$ que nous appelons \textit{tresse normale}. 

\vspace{0.2cm}

Nous d\'emontrons ensuite (Th\'eor\`eme \ref{decollement}) que si l'on a un 
plongement lisse $\pi:\overline{\Sigma}\rightarrow V$, 
alors il existe un plongement lisse 
$\pi_{\infty} :\mathcal S_g\rightarrow V$ du sol\'eno\"{\i}de de Sullivan de genre $g=g(\Sigma)$ 
aussi proche que l'on veut de $\pi$ dans la topologie lisse. 

\vspace{0.2cm}

Bien que la notion de surface branch\'ee holomorphe n'ait pas de sens \`a cause du principe de prolongement analytique,
celle de surface branch\'ee $\varepsilon$-holomorphe en a lorsque $\varepsilon$ est strictement positif. 
Pour tout $\varepsilon>0$, nous construisons un plongement lisse $\varepsilon$-holomorphe de 
$\overline{\Sigma}_2$ dans le plan projectif complexe (Th\'eor\`eme \ref{plgt surface branchee}).
On en d\'eduit que pour tout $\varepsilon>$, 
il existe un plongement lisse $\varepsilon$-holomorphe du sol\'eno\"{\i}de de Sullivan de genre $2$ dans ${\bf C}P^2$
(Corollaire \ref{plgt solenoide}).\\

\textit{Remerciements.} Je remercie \'Etienne Ghys et 
Andr\'e Haefliger pour l'encouragement qu'ils m'ont apport\'e. Je remercie aussi Sidney Frankel, Emmanuel Giroux, 
Alexei Glutsyuk, Martin Pinsonnault, Jean-Claude Sikorav et Jean-Yves Welschinger 
pour les nombreuses conversations que l'on a \'echang\'ees sur ce sujet. 
Ce travail a \'et\'e rendu possible gr\^ace \`a l'hospitalit\'e de l'universit\'e de Gen\`eve et de l'universit\'e de Toronto. 

\section{Sol\'eno\"{\i}des et surfaces branch\'ees}\label{solenoide de Sullivan}

Une \textit{lamination par surfaces de Riemann} d'un espace topologique $X$ est un atlas $\mathcal L$ constitu\'e
d'hom\'eomorphismes $\varphi :U \rightarrow {\bf D} \times T$ (les ouverts $U$ forment un recouvrement de $X$) \`a valeurs dans le produit du disque 
unit\'e par un espace m\'etrique, tel que:

\ (i) Les changements de cartes pr\'eservent la fibration locale par disques. 

(ii) Ils sont holomorphes le long des fibres.

\noindent Lorsque les espaces transverses $T$ sont totalement discontinus, nous appelons la lamination un \textit{soleno\"{\i}de}. 

\subsection{Le sol\'eno\"{\i}de de Sullivan}
Soit ${\bf D}^*=\{ 0<|z|<1\}$ le disque \'epoint\'e et $r$ le rev\^etement double \[r:z\in {\bf D}^* \mapsto z^2 \in {\bf D}^*.\]  
Appelons \textit{orbite} une partie $\{z_n\}_{n\in {\bf Z}}$ de $\bf D^*$ telle que pour tout entier $n$, 
\[ z_{n+1}= z_n ^2.\] 
telle que pour tout point $z$ de $O$ les propri\'et\'es suivantes sont v\'erifi\'ees. 
Les orbites de $r$ sont des parties discr\`etes de ${\bf D}^*$, 
puisque les it\'er\'es positifs de $r$ convergent vers $0$, et les branches inverses 
de $r$ convergent vers le cercle unit\'e.
L'ensemble de toutes les orbites de $r$, muni de la topologie de Hausdorff 
sur les compacts, est un espace topologique compact que nous notons 
$\mathcal S$. 

\begin{proposition}\label{structure solenoidale} 
L'espace $\mathcal S$ a une structure de sol\'eno\"{\i}de par surfaces de Riemann. \end{proposition}

\textit{D\'emonstration.} Soit $z$ un point du disque \'epoint\'e. Une orbite $O$ de $r$ contenant $z$ 
est construite de la fa\c{c}on suivante. 
On choisit une d\'etermination $z_1$ de la racine de $z$, une d\'etermination $z_2$ de la racine de $z_1$ etc. L'orbite $O$ est la r\'eunion 
des it\'er\'es de $z$ par $r$ et de la famille $\{z_1, z_2, \ldots\}$. Ainsi, l'ensemble des orbites de $r$ qui contiennent le point $z$ 
s'identifie au bord de l'arbre diadique dessin\'e sur la Figure~\ref{arbre}. 
Cet espace est donc hom\'eomorphe \`a l'ensemble de Cantor.

\begin{figure}[h!]\label{arbre}
\begin{center} 
\input{arbre.pstex_t}
\caption{}
\end{center}
\end{figure}

Consid\'erons un disque $\delta $ de ${\bf D}^*$, assez petit pour que toute orbite de $r$ l'intersecte en 
au plus un point, et analysons la structure de l'ouvert $B(\delta)$ 
des orbites $O$ de $r$ qui intersectent $\delta$. 
Notons $z(O)$ le point d'intersection d'une telle orbite $O$ avec $\delta$. L'application 
\[  z: O\in B(\delta)  \mapsto z(O) \in \delta \]
est continue, et ses fibres sont des ensembles de Cantor. De plus, \'etant donn\'e un point $O$ de $B(\delta)$, 
nous pouvons suivre les d\'eterminations
analytiques des branches inverses de $r$ lorsqu'on fait varier $z$ dans $\delta$. 
Ceci prouve que $z$ est une fibration localement triviale. 
Il existe donc une application $t: B(\delta )\rightarrow T$, o\`u $T$ est l'ensemble de Cantor,
telle que $(z,t)$ est un hom\'eomorphisme de $B(\delta)$ dans $\delta \times T$. 
Ce sont les cartes qui munissent $\mathcal S$ d'une structure 
de sol\'eno\"{\i}de par surfaces de Riemann. La Proposition \ref{structure solenoidale} est d\'emontr\'ee. \\

\begin{figure}[h!]\label{surfacesigmainfty}
\begin{center}
\input{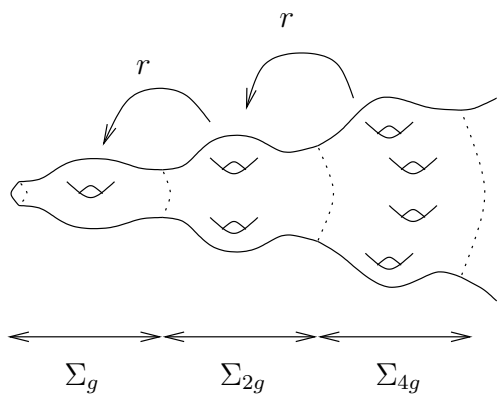}
\caption{$\Sigma_{\infty,g}$}
\end{center}
\end{figure}

Faisons une construction analogue en genre non nul. 
Consid\'erons une suite de surfaces \textit{hyperboliques} compactes $\Sigma_{2^ng}$ de genre $2^n g$, dont le bord
est form\'e de deux composantes connexes g\'eod\'esiques $\partial_- \Sigma_{2^ng}$ et $\partial _+\Sigma_{2^ng}$.
Et des rev\^etements doubles localement isom\'etriques $r_n: \Sigma_{2^{n+1}g} \rightarrow \Sigma_{2^ng}$. 
Formons une surface hyperbolique non compacte $\Sigma_{\infty,g}$ en attachant les surfaces $\Sigma_{2^ng}$ 
par une identification g\'eod\'esique
de $\partial _-\Sigma_{2^{n+1}g}$ avec $\partial _+ \Sigma _{2^ng}$ (voir Figure~\ref{surfacesigmainfty}). 
Il y a plusieurs fa\c{c}ons de coller ces surfaces, mais nous faisons en sorte de 
construire aussi un rev\^etement localement isom\'etrique 
\[ r: \Sigma_{\infty,g}-\Sigma_g \rightarrow \Sigma _{\infty,g},\] 
qui en restriction \`a $\Sigma_{2^{n+1}g}$ est le rev\^etement 
double $r_n$ au dessus de $\Sigma_{2^ng}$.

Une \textit{orbite} de $r$ est une partie $O$ de $\Sigma_{\infty,g}$ telle que pour tout point $z$ de $O$, 
les deux propri\'et\'es suivantes soient satisfaites. 
D'une part, les it\'er\'es par $r$ de $z$, lorsqu'ils sont d\'efinis, sont dans $O$, et d'autre part, 
il existe un unique $w$ dans $O$ tel que $r(w)=z$.   
L'espace de toutes les orbites de $r$, muni de la topologie de Hausdorff sur les compacts, 
est un sol\'eno\"{\i}de par surfaces de Riemann 
que l'on note $\mathcal S_g$. La d\'emonstration de ce fait est analogue \`a celle de la Proposition~\ref{structure solenoidale}.   
D'ailleurs, lorsque $g$ est nul, le sol\'eno\"{\i}de $\mathcal S_0$ est biholomorphe au sol\'eno\"{\i}de de Sullivan $\mathcal S$.

\subsection{Approximation de $\mathcal S_g$ par des surfaces branch\'ees}

Soit $g$ un entier positif. Pour tout entier $n$, le rev\^etement $r$ induit un rev\^etement double  
\[    r: \partial_+ \Sigma_{2^ng} \rightarrow \partial_- \Sigma_{2^ng} .\]
Notons 
\begin{equation}\label{involution} i:\partial_+ \Sigma_{2^ng} \rightarrow \partial_+ \Sigma_{2^ng}\end{equation} 
l'involution d\'efinie par ce rev\^etement, et 
identifions les points $z$, $i(z)$ et $r(z)$ de $\Sigma_{2^ng}$ si $z$ appartient \`a $\partial_+ \Sigma_{2^ng} $. 

\begin{figure}[h!]\label{overlinesigma}
\begin{center}
\input{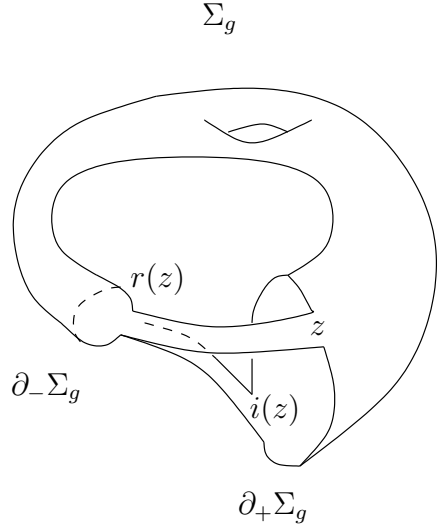}
\caption{$\overline{\Sigma_g}$}
\end{center}
\end{figure}

Nous obtenons une \textit{surface branch\'ee} $\overline{\Sigma_{2^ng}}$, dont le lieu de branchement est un cercle 
localement hom\'eomorphe \`a un livre ouvert \`a trois feuilles. 
En effet, le voisinage d'un point du lieu de branchement est obtenu en attachant 
trois demi-disques ouverts $D(r(z))$, $D(z)$ et 
$D(i(z))$ (voir Figures \ref{overlinesigma} et \ref{branchement}). Lorsqu'on fait le tour du lieu de branchement, 
les deux feuillets $D(z)$ et $D(i(z))$ sont \'echang\'es alors que $D(r(z))$ reste fixe.

Chaque orbite $O$ de $r$ intersecte $\Sigma_{2^ng}$ en un point \`a l'int\'erieur, ou en deux points $r(z)$ et $z$ sur chaque bord de 
$\Sigma_{2^ng}$. L'application 
\begin{equation} \label{applications pgn} p_{g,n} :O\in {\mathcal S_g} \mapsto O\cap 
\Sigma_{2^n g}\in \overline{\Sigma_{2^ng}} \end{equation}
est donc bien d\'efinie et continue. 
\begin{lemme} \label{approximation} 
Le diam\`etre des fibres de $p_{g,n}$ tend uniform\'ement vers $0$ lorsque $n$ tend vers 
l'infini.\end{lemme}

\begin{remarque} L'espace topologique $\mathcal S_g$ est m\'etrisable. Le r\'esultat est 
ind\'ependant de la m\'etrique que nous choisissons. \end{remarque}

\textit{D\'emonstration.} Soit $n$ un entier positif. L'image $p_n(O)$ d'une orbite de $O$ d\'etermi\-ne 
son intersection avec la r\'eunion 
\[\bigcup _{0\leq k\leq n-1}  \Sigma_{2^kg} \subset \Sigma_{\infty,g}.\] 
Cette famille de parties de $\Sigma_{\infty,g}$ est croissante et sa r\'eunion est $\Sigma_{\infty,g}$,
donc le diam\`etre des fibres de l'application 
$p_n$ tend uniform\'ement vers $0$ lorsque $n$ tend vers l'infini.

\subsection{Structure lisse}\label{structure lisse}
Soit $\mathcal L$ une lamination par surfaces de Riemann d'un espace topologique $X$. 
Une fonction $f:X\rightarrow {\bf R}$ est de classe $C^k$ si dans une carte $(z,t):U\rightarrow {\bf D}\times T$, 
les fonctions 
\[    f_t(z)= f(t,z) : {\bf D}\rightarrow {\bf R}\]
sont de classe $C^k$ et d\'ependent continument de $t$ dans la topologie $C^k$.  
Cette d\'efinition est bien ind\'ependante du syst\`eme de coordonn\'ees choisi dans $\mathcal L$, car les changements 
de cartes sont holomorphes.

\vspace{0.2cm}

Il s'agit de construire une structure lisse sur les surfaces branch\'ees $\overline{\Sigma_{2^ng}}$, qui approche la structure 
lisse de $\mathcal S_g$ lorsque $n$ tend vers l'infini
(voir Lemme \ref{approximation lisse}). Cette structure est d\'efinie par 
l'ensemble de ses fonctions lisses \`a valeurs dans ${\bf R}$. Elle v\'erifie 
la propri\'et\'e suivante: si $f:\overline{\Sigma_{2^ng}}\rightarrow 
{\bf R}^q$ et $g :{\bf R}^q\rightarrow {\bf R}$ sont des fonctions lisses, alors $g\circ f$ est une fonction lisse.

\begin{definition} La structure lisse au voisinage du lieu de branchement est obtenue de la fa\c{c}on suivante. 
Il existe un voisinage $\overline{A}$ du lieu de 
branchement de $\overline{\Sigma_{2^ng}}$ qui est hom\'eomorphe au quotient de l'anneau 
$A=[0,1]\times {\bf S}^1$ par la relation qui identifie 
les points $(r,z)$ et $(r,-z)$ si $r$ est sup\'erieur \`a $1/2$. Notons 
\[ p: A\rightarrow {\overline A}\]
la projection canonique. Une fontion $f: \overline{A} \rightarrow {\bf R}$ est dite de classe $C^k$ si $f\circ p$ est de classe $C^k$
($k=1,\ldots,\infty$).

En dehors du lieu de branchement, la structure lisse est l'unique structure lisse sur une surface topologique \`a bord. \end{definition}

\begin{figure}[h!]\label{branchement}
\begin{center}
\input{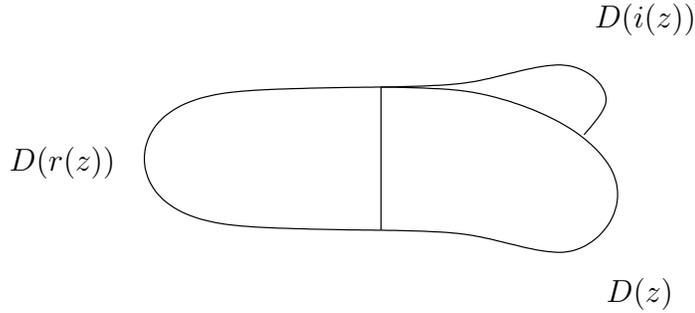}
\caption{Lieu de branchement}
\end{center}
\end{figure}

Notons qu'avec cette d\'efinition, il est possible de d\'efinir un espace tangent en tout point de $\overline{\Sigma_{2^ng}}$, y compris 
sur le lieu de branchement. La Figure \ref{branchement} illustre ce fait~: elle repr\'esente 
l'image d'un voisinage du lieu de branchement par une immersion lisse \`a valeurs dans ${\bf R}^3$.  

Pour tous les entiers $g$ et $n$, les applications $p_{g,n} :{\mathcal S_g}\rightarrow \overline{\Sigma_{2^ng} }$ 
d\'efinies en~\ref{applications pgn} sont lisses. 
En fait, une fonction $f: \overline{\Sigma_{2^n g} } \rightarrow {\bf R}$ 
est lisse si et seulement si $f\circ p_{g,n}$ est une fonction lisse. 

\begin{lemme}\label{approximation lisse} Pour toute fonction lisse $f:{ \mathcal S_g}\rightarrow {\bf R}$,  
il existe une suite de fonctions lisses $f_n : \overline{\Sigma_{2^n g}} \rightarrow {\bf R}$
telles que $f_n \circ p_{g,n}$ convergent vers $f$ dans la topologie $C^{\infty}$ lorsque $n$ 
tend vers l'infini.\end{lemme}

\textit{D\'emonstration.} Il suffit de d\'emontrer que, pour tout entier $k$ positif, 
il existe une suite de fonctions lisses $f_n$ telles que 
$f_n\circ p_{g,n}$ convergent vers $f$ dans la topologie $C^k$. 

Nous affirmons qu'il est possible de trouver un recouvrement ouvert fini 
$\{U_k\}$ de $\mathcal S_g$ et des cartes holomorphes $U_k \rightarrow {\bf D}\times T_k$ telles que 
les fibres de $p_{g,0}$ soient contenues dans les ensembles $\{z\}\times T_k$. Ceci d\'ecoule directement 
de la d\'emonstration de~\ref{structure solenoidale}. Observons alors que
pour tout $n$ et tout $k$, nous pouvons trouver une partition (finie) de $T_k$ en des parties ouvertes et ferm\'ees 
$T_{k,l}^n$ dont les diam\`etres tendent uniform\'ement vers $0$ lorsque $n$ tend vers l'infini, et
telles que les fibres de $p_{g,n}$ soient contenues dans les ensembles $z\times T_{k,l}^n$.

Soit $0<\rho <1$ un r\'eel
tel que, en notant ${\bf D}_{\rho} =\{ |z|\leq \rho\}$, les ouverts ${\bf D}_{\rho}\times T_k$ recouvrent encore $\mathcal S_g$. 
Quitte \`a prendre une partition de l'unit\'e form\'ee de fonctions dont les supports sont inclus dans ${\bf D}_{\rho}\times T_k$,
nous pouvons supposer que le support de $f$ est contenu dans l'un des ${\bf D}_{\rho}\times T_k$. Pour tout $n$ et 
tout $l$, choisissons un point $t_{k, l}^n$ dans $T_{k,l}^n$, et posons 
\[  g_n(z,t) = f(z,t_{k,l}^n) \]
si $t$ appartient \`a $T_{k,l}^n$, et $f_n(z,t)=0$ sinon. La fonction $g_n$ est une fonction lisse. De plus, puisque 
le diam\`etre de $T_{k,l}^n$ tend uniform\'ement vers $0$ lorsque $n$ tend vers l'infini, les fonctions $g_n$ 
convergent vers $f$ dans la topologie $C^k$ lorsque $n$ tend vers l'infini. 

Les fibres de $p_n$ sont contenues 
dans les ensembles $z\times T_{k,l}$. Il existe donc une fonction 
$f_n:\overline{\Sigma_{2^n g}}\rightarrow {\bf R}$ telle que $g_n=f_n\circ p_n$. Le Lemme \ref{approximation lisse} est d\'emontr\'e.

\section{Voisinage tubulaire des surfaces branch\'ees $\overline{\Sigma_g}$}\label{voisinage}

Soit $g$ un entier positif. 
Dans ce paragraphe nous \'etudions l'image $\overline{\Sigma}=\pi (\overline{\Sigma_g})$ d'un plongement lisse 
$\pi :\overline{\Sigma}_g\rightarrow V$ dans une 
vari\'et\'e \textit{orient\'ee} $V$ de dimension $4$. Nous appelons \textit{plongement lisse} un plongement topologique dont la d\'eriv\'ee 
est injective, y compris sur le lieu de branchement.

\subsection{Tresse normale}
Un voisinage d'une sous-vari\'et\'e compacte est, d'apr\`es le Th\'eor\`eme du voisinage tubulaire, diff\'eomorphe 
\`a son fibr\'e normal. Le cas des surfaces branch\'ees est tout \`a fait diff\'erent~: 
le fibr\'e normal de $\overline{\Sigma}$ est toujours trivial, mais pourtant 
il existe un autre invariant, que nous appelons \textit{tresse normale}. 
Cette terminologie sera justifi\'ee en \ref{terminologie}. 

\begin{proposition}\label{tout fibre est trivial} Le fibr\'e normal de $\overline{\Sigma} \subset V$  
est trivial. En fait, $H^2(\overline{\Sigma}, {\bf Z})=0.$\end{proposition}

\textit{D\'emonstration.} Soit $\overline{A}$ un voisinage du lieu de branchement de $\overline{\Sigma}$, et $\Sigma$ l'ext\'erieur du 
lieu de branchement. Nous avons $\overline{\Sigma} = \overline{A} \cup \Sigma$, et la suite exacte de Mayer-Vietoris s'\'ecrit~:
\[  H^1(\Sigma) \oplus H^1(\overline{A}) \stackrel{i_{\Sigma}^*\oplus i_{\overline{A}}^*}{\longrightarrow} H^1(\Sigma\cap \overline{A})
\rightarrow H^2(\overline{\Sigma}) 
\rightarrow H^2(\Sigma)\oplus H^2(\overline{A}).\]  
Remarquons d\'ej\`a que le dernier terme de cette suite est nul, d'une part parce que le $H^2$ d'une surface non compacte est nul,
et d'autre part parce que $\overline{A}$ est homotopiquement \'equivalent \`a un cercle. 

Il nous suffit alors de montrer que $i_{\Sigma}^*\oplus i_{\overline{A}}^*$ est surjective. 
Soit $c_+$ et $c_-$ des courbes obtenues en poussant $\partial_- \Sigma$ et $\partial_+\Sigma$ \`a 
l'int\'erieur de $\Sigma$. Comme $\Sigma \cup \overline{A}$ est homotopiquement \'equivalent \`a la r\'eunion 
$c_-\cup c_+$, le groupe $H^1(\Sigma \cup \overline{A} )$ est le groupe ab\'elien libre engendr\'e par $[c_- ]^*$ et $[c_+]^*$
($[c_- ]^*,[c_+]^*$ est la base duale de $[c_- ],[c_+]$). Mais comme les courbes $c_-$ et $c_+$ sont homotopes 
aux deux bords de $\Sigma$, il existe une forme lin\'eaire $\lambda _{\Sigma}$ sur $H_1(\Sigma)$ qui vaut $1$ sur $[c_+]$ et $1$ ou $-1$ sur $[c_-]$.
D'autre part, $H_1(\overline{A})= {\bf Z} [c_-]$, et $[c_+]= 2[c_-]$~; il existe donc une forme lin\'eaire sur $\lambda _{\overline{A}}$ 
sur $H_1(\overline{A})$ qui vaut $1$ sur $[c_-]$ et $2$ sur $[c_+]$. Les formes $i_{\Sigma}^* \lambda_{\Sigma}$ et $i_{\overline{A}}^*\lambda_{\overline{A}}$
engendrent donc $H^1(\Sigma\cap \overline{A})$ et la Proposition est d\'emontr\'ee. \\ 
  
\begin{figure}[h!] \label{sectionnormale}
\begin{center}
\input{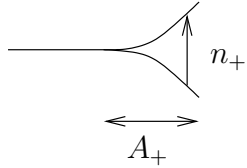}
\caption{La section $n_+$}
\end{center}
\end{figure}

D\'efinissons \`a pr\'esent la tresse normale.  
D'apr\`es la Proposition \ref{tout fibre est trivial}, 
il existe une section continue $n:\overline{\Sigma} \rightarrow N$
du fibr\'e normal de $\overline{\Sigma}$, qui ne s'annule pas. 
Comme $c_+$ est homologue \`a $0$ dans $\overline{\Sigma}$, la restriction $n_{|c_+} $ de $n$ \`a $c_+$ est uniquement d\'etermin\'ee 
\`a homotopie pr\`es. 

Nous avons d'autre part une section $n_+: A_+ \rightarrow N$ d\'efinie par
\[ \pi\circ i(x)= \exp (n_+ (i(x)) )  ,\]
et que nous avons repr\'esent\'ee sur la Figure \ref{sectionnormale}. Elle mesure comment les feuillets tournent au voisinage du lieu de branchement 
de $\overline{\Sigma}$. Rappelons que l'application $i$ est une involution sans point fixe d\'efinie sur $A_+$ (voir \ref{involution}). L'application 
$\exp$ est l'application exponentielle d\'efinie sur le fibr\'e normal de $\overline{\Sigma}$ (identifi\'e \`a l'orthogonal du fibr\'e tangent) 
\`a valeurs dans $V$. La section $n_+$ ne s'annule pas et sa classe d'homotopie est aussi bien d\'efinie.

\begin{definition}
La \textit{tresse normale} $N(\overline{\Sigma},V) $ est l'indice de la fonction 
\[\frac{n_+}{n}: c_+\rightarrow {\bf C}^*,\]
l'orientation de $c_+$ \'etant donn\'ee par celle de $\partial _-\Sigma$ (ou l'orientation oppos\'ee \`a $\partial_+\Sigma$). \end{definition}

\begin{remarque}\label{terminologie} La tresse normale peut \^etre beaucoup plus compliqu\'ee pour 
des surfaces branch\'ees dont le lieu de ramification est un livre ouvert \`a 
plus de quatre feuilles. Consid\'e\-rons 
par exemple la surface branch\'ee obtenue en quotientant $[0,1]\times {\bf S}^1$ par la relation qui identifie 
les points $(r,z)$, $(r,jz)$ et $(r,j^2z)$ lorsque $r\leq 1/2$ (ici $j=\exp(2i\pi/3)$). 
Elle est munie de la structure lisse naturelle pour laquelle les fonctions lisses sont les fonctions 
lisses sur $[0,1]\times {\bf S}^1$ qui sont constantes sur les classes d'\'equivalence. 
Le bord de cette surface est form\'ee de deux cercles que nous attachons par un diff\'eomorphisme, 
et nous obtenons de la sorte une surface branch\'ee $\overline{\Sigma}^{(3)}$. 
Le lecteur pourra s'assurer qu'il est possible de d\'efinir un invariant 
\[ N: \{ \mathrm {plongements\ lisses}\ \pi:\overline{\Sigma}^{(3)}\rightarrow V\}\rightarrow B_3\]
\`a valeurs dans le groupe des tresses \`a $3$ brins $B_3$. C'est la raison pour laquelle 
nous avons appel\'e l'invariant $N$ la tresse normale. 

Un corollaire du Th\'eor\`eme \ref{voisinage global} est le fait que la tresse normale prend des valeurs impaires pour les surfaces 
$\Sigma_m^{(2)}$. Il n'y a rien de surprenant \`a cela,
puisque $B_2$ est isomorphe \`a ${\bf Z}$, 
les tresses impaires correspondant aux tresses \textit{connexes}.\end{remarque} 

\subsection{Formes normales}\label{forme normale} Pour tout entier $p$, nous construisons 
une vari\'et\'e orient\'ee $V_p$ de dimension $4$ 
et une surface branch\'ee lisse orient\'ee $\overline{\Sigma}\subset V_p$ diff\'eomor\-phe \`a $\overline{\Sigma_g}$ de tresse normale 
$N(\overline{\Sigma}, V_p) =1-2p$. Ceci montre que la tresse normale n'est pas un invariant trivial.  

 La vari\'et\'e $V_p$ a du bord, des coins et n'est pas compacte. 
C'est la r\'eunion des deux vari\'et\'es $W_1= A\times {\bf C}$ ($A=[0,1]\times {\bf S}^1$) et $W_2= \Sigma_g \times {\bf D}$
que l'on a coll\'ees sur des parties de leur bord par les plongements lisses~:
\[ \Phi_+:(z,\zeta)\in \partial _+\Sigma_g \times {\bf D} \mapsto (z^2,4z+z^p\zeta) \in \partial _+ A\times {\bf C},\]
et 
\[ \Phi_- : (z,\zeta)\in \partial_-\Sigma_g \times {\bf D}\mapsto (z,\zeta) \in \partial_-A\times {\bf C},\]   
les coordonn\'ees $z$ en restriction aux bords de $A$ et $\Sigma$ pr\'eservant l'orientation de 
$\partial_-\Sigma,\ \partial_+A$, et renversant l'orientation de $\partial_+\Sigma,\ \partial_-A$. 

\begin{figure}[h!]\label{formenormale} 
\begin{center}
\input{formenormale.pstex_t}
\caption{$V_p$}
\end{center}
\end{figure}

La surface branch\'ee $\overline{\Sigma}\subset V_p$ est construite de la mani\`ere suivante. 
Soient \[A_- = [0,1/2)\times {\bf S}^1\ \ \ \ \ \ \ \ A_+= (1/2,1]\times {\bf S}^1,\]
et $f: A\rightarrow {\bf C} $ une fonction lisse v\'erifiant les propri\'et\'es suivantes~:

\ \ (i) elle est nulle sur $A_-$, 

\ (ii) elle ne s'annule pas sur $A_+$, 

(iii) elle vaut $16z$ sur un voisinage de $\{1\}\times {\bf S}^1$,

(iv) au voisinage de $\{1/2\}\times {\bf S}^1$, les d\'eterminations de sa racine carr\'ee sont lisses.

Soient $\overline{A}\subset W_1$ la surface branch\'ee d\'efinie par l'\'equation 
\[  \overline{A} =\{  \zeta^2 =f \},\]
et $\Sigma = \Sigma_g\times 0 \subset W_2$.
La r\'eunion $\overline{\Sigma}= \Sigma \cup \overline{A} \subset V_p$ est l'image d'un plongement lisse de
$\overline{\Sigma_g}$. 

\begin{remarque}\label{non rigidite}
Il n'est pas vrai que si $U= \overline{D(z)}\cup\overline{D(i(z))}\cup \overline{D(r(z))}$,
et si $\pi_i: U\rightarrow {\bf R}^q, \ i=0,1$ sont deux plongements lisses, 
alors il existe un diff\'eomorphisme $\phi: {\bf R}^q\rightarrow {\bf R}^q$
tel que $\phi(\pi_0(U)) = \pi_1(U)$. 
Pour construire un contre-exemple $1$-dimensionel, d\'efinissons
\[  B_{\psi} = \{(x,y)| \ y^2=\psi(x)\}\]
pour toute fonction $\psi: {\bf R}\rightarrow {\bf R}$. 
Nous affirmons qu'il n'existe 
pas de diff\'eomorphisme de classe $C^1$ de ${\bf R}^2$ dans lui-m\^eme qui envoit 
$B_{\psi_0}$ sur $B_{\psi_1}$, si $\psi_0$ et $\psi_1$ sont choisies de la fa\c{c}on suivante~: $\psi_0$ et $\psi_1$ 
sont identiquement nulle sur $(-\infty,0]$, et pour $x>0$~:
\[  \psi_0(x)=\exp(-1/x),\ \ \ \ \ \  \psi_2(x)= (2+\sin(1/x)) \exp(-1/x).\]
Pourtant, $B_{\psi_i},\ i=0,1$ sont des images d'un plongement lisse du graphe branch\'e lisse $\mathcal G$, obtenu en identifiant 
deux points $(0,t)\in 0\times [0,1]$ et $(1,t)\in [0,1]$ si $t\leq 1/2$. Une fonction lisse sur $\mathcal G$ 
\'etant une fonction qui se rel\`eve en une fonction lisse sur chaque $i\times [0,1],\ i=0,1$.  
Ainsi, les formes normales que nous avons propos\'ees ne sont pas toutes diff\'eomor\-phes, m\^eme \`a $g$ et $p$ 
fix\'e~: leur classe de diff\'eomorphie d\'epend de la fonction $f$. 
\end{remarque}

\begin{proposition} La tresse normale $N(\overline{\Sigma}, V_p)$ est \'egale \`a $1-2p$. \end{proposition}

\textit{D\'emonstration.} Commen\c{c}ons par calculer une section du fibr\'e normal de $\overline{\Sigma}$ qui ne 
s'annule pas. Consid\'erons des fonctions lisses $z:\Sigma_g\rightarrow {\bf S}^1$ et $z:A\rightarrow {\bf S}^1$ qui prolongent les 
coordonn\'ees $z:\partial .\rightarrow {\bf S}^1$ d\'efinies sur les bords de $\Sigma$ et $A$. 
D\'efinissons les sections $n_1:\overline{A}\rightarrow N$ et $n_2:\Sigma\rightarrow N$ par  
\[  n_1  = z^{\alpha} \frac{\partial} {\partial \zeta}\ \ \ \mathrm{et}\ \ \ n_2 =  z^{\beta} \frac{\partial }{\partial \zeta},\]
o\`u $\alpha $ et $\beta$ sont des entiers. Nous avons 
\[d\Phi_+.n_2 (z^2,z) = z^{\beta+p} \frac{\partial}{\partial \zeta} \ \ \ \mathrm{et}\ \ \ d\Phi_- .n_2 (z,0)=z^{\beta}\frac{\partial}{\partial \zeta}.\] 
Pour que $n_1$ et $n_2$ d\'efinissent une section globale $n$, nous devons choisir $\alpha=\beta=p$. 

La section $n_+$ est donn\'ee sur $c_+$ par 
\[ n_+(z^2,z) = -2z \frac{\partial}{\partial \zeta},\]
et nous avons donc pour tout point $(z,z^2)$ de $c_+$,  
\[ n_+(z^2,z)=-2z^{1-2p}n(z^2,z).\]
La Proposition est donc d\'emontr\'ee.

\subsection{Topologie du voisinage de $\overline{\Sigma}$}
Nous montrons maintenant que la topologie d'un voisinage de l'image $\overline{\Sigma}=\pi (\overline{\Sigma_g})$
d'un plongement lisse $\pi: \overline{\Sigma _g}\rightarrow V$ est compl\`etement d\'etermin\'ee par sa tresse normale, et que 
$\overline{\Sigma}$ est d\'ecrite par l'un des mod\`eles que nous avons 
propos\'e~(donn\'e par le choix d'une fonction $f:A\rightarrow {\bf C}$).  
\begin{theoreme}\label{voisinage global} Soit $V^4$ une vari\'et\'e compacte de dimension $4$ orient\'ee 
et $\pi : \overline{\Sigma_g} \rightarrow V$ un plongement lisse. Il existe un entier $p$, un voisinage $\mathcal V$ 
de $\overline{\Sigma}=\pi(\overline{\Sigma_g})$ et un diff\'eomorphi\-sme $\Phi : \mathcal V\rightarrow V_p$ pr\'eservant l'orientation 
tel que $\Phi (\overline{\Sigma})$ est l'une des surfaces branch\'ees construites en \ref{forme normale}.\end{theoreme}

\textit{D\'emonstration.} Nous commen\c{c}ons par \'etudier le voisinage du lieu de branchement. 

\vspace{0.2cm}

{\bf Fait I.} \textit{Il existe un voisinage $\mathcal W_1$ du lieu de branchement de $\overline{\Sigma}$ et des coordonn\'ees 
$\psi_1=(x,\zeta) : \mathcal W_1 \rightarrow A\times \bf C$ tels que 
\[ \overline{\Sigma} \cap {\mathcal W}_1 = \{  (x,\zeta) |\ \zeta^2 = f(x) \},\]
o\`u $f: A\rightarrow \bf C$ est une fonction lisse v\'erifiant les conditions (i), (ii), (iii) et (iv).} 

\vspace{0.2cm}

\textit{D\'emonstration.} Soit $b\subset \overline{\Sigma}$ le lieu de branchement de $\overline{\Sigma}$. Comme le fibr\'e tangent de 
$\overline{\Sigma} $ est orient\'e, 
il existe un anneau lisse plong\'e $A\subset V$ tel que $T_b A = T_b \overline{\Sigma}$, et s\'epar\'e par $b$ en deux anneaux $A_+$ et $A_-$. 
Puisque $V$ est orient\'ee, $A$ est aussi 
transversalement orient\'e, et en vertu du Th\'eor\`eme du voisinage tubulaire, il existe un voisinage $\mathcal W_1$ 
de $b$ et un diff\'eomorphisme $(x,\zeta') : \mathcal W_1 \rightarrow A\times \bf C$ tel que $\psi _{|A} = id \times 0$. 
Lorsqu'on choisit l'anneau $A$ assez fin, et quitte \`a permuter les anneaux $A_+$ et $A_-$, 
l'intersection de $\overline{\Sigma} $ avec une fibre $x\times \bf C$ est form\'ee de deux points 
$(x, \zeta'_i(x))$, $i=1,2$, \'egaux si $x\in A_-$ et distincts si $x\in A_+$. Au voisinage d'un point $x$ de $A$, ces points sont param\'etr\'es
par des d\'eterminations 
\[ \zeta'_i: A\rightarrow {\bf C},\ \ \ \ i=1,2\]
qui sont des fonctions lisses. De plus la topologie du lieu de branchement impose que $\zeta' _1$ soit chang\'ee en $\zeta'_2$ 
lorsqu'on fait le tour de $A_+$. La fonction 
\[  f' = (\frac{\zeta'_1 -\zeta'_2}{2})^2\]
est donc une fonction lisse d'indice impair $2m+1$. Soit $z:A \rightarrow {\bf C}^*$ une fonction d'indice $1$. 
Consid\'erons les coordonn\'ees 
\[  (x,\zeta) = (x, z^{-m} (\zeta' -\frac{\zeta'_1(x)+\zeta'_2(x)}{2})).\]
Dans ces coordonn\'ees, l'\'equation de $\overline{\Sigma}$ est donn\'ee par 
\[  \{  \zeta^2 = f(x) \},\] 
o\`u $f= z^{-2m} f'$, qui est bien d'indice $1$ en restriction \`a $A_+$. Le Fait I est d\'emontr\'e.

\vspace{0.2cm}

La surface $\overline{\Sigma} - \mathrm{Int}(\mathcal W_1)$ est diff\'eomorphe \`a $\Sigma_g$, et comme $V$ est orient\'ee, 
son fibr\'e normal est trivial. Soit $\partial _+ \Sigma = \{ \zeta ^2=16 z\}$, et $\partial _- \Sigma = \partial_- A\times 0$
dans les coordonn\'ees du Fait I.  
Observons que nous avons des trivialisations $s_-:\partial_-\Sigma_g\rightarrow N$ et $s_+:\partial_+\Sigma_g\rightarrow N$, 
donn\'ees par le champ de vecteurs $ \frac {\partial}{\partial \zeta}$. 

\vspace{0.2cm}

{\bf Fait II.} \textit{Il existe une section lisse $s:\Sigma_g \rightarrow N$ et un entier $p$ tels que $s_{|\partial_- \Sigma_g}=s_-$ et 
\[  s(0,z,\zeta) =  z^p s_+(0,z,\zeta),\]
pour tout $(0,z,\zeta)$ de $\partial _+ \Sigma_g$.}

\vspace{0.2cm}

\textit{D\'emonstration.} Le fibr\'e normal $N\rightarrow \Sigma_g$ est trivial, et la courbe 
$\partial _-\Sigma_g$ est non homologue \`a 
$0$, donc il existe une section lisse $\overline{s}:\Sigma_g\rightarrow N$ qui prolonge $s_-$ et qui ne s'annule pas. 
Prolongeons $s_+$ en une section 
lisse sur un anneau $C_+$ bordant $\partial_+\Sigma_g$. Soit $f= \overline{s}/s_+$ et $p$ le degr\'e de $f$. 
Comme $f^{-1} z^{p}$ est de degr\'e nul, elle est homotope \`a la fonction constante $1$. Il existe donc une fonction
lisse $h:\Sigma_g\rightarrow {\bf C}^*$ 
qui vaut $f^{-1}z^{p}$ sur $\partial _+\Sigma_g$ et $1$ \`a l'ext\'erieur de $C_+$. 
La section $s=h\overline{s}$ v\'erifie alors le Fait II.

\vspace{0.2cm}

Achevons la d\'emonstration du Th\'eor\`eme \ref{voisinage global}. 
Consid\'erons une m\'etri\-que riemannienne sur $V$ qui est la m\'etrique riemannienne 
plate $\psi_1^*(dr^2+|dz|^2+|d\zeta|^2)$ en res\-triction \`a $\mathcal W_1$. 
Nous pouvons supposer que $\overline{A}$ est orthogonal au bord $\partial_+A \times {\bf C}$. 
Identifions le fibr\'e normal de $\overline{\Sigma}$ avec 
l'orthogonal de son fibr\'e tangent. Lorsque $\alpha>0$ est assez petit, l'application 
\[ \psi_2 : (x,\zeta) \in \Sigma \times {\bf D} \mapsto \exp (\alpha \zeta s(x))\in V\]
est bien d\'efinie et c'est un diff\'eomorphisme sur son image $\mathcal W_2$. 
D'autre part, $\mathcal W_1$ et $\mathcal W_2$ s'intersectent exactement sur $\psi_2(\partial \Sigma \times \overline {\bf D}^2)$. 
Il suffit alors d'\'etudier les changements de cartes $\psi_2\circ\psi_1^{-1}$ pour achever la d\'emonstration du Th\'eor\`eme.

\section{D\'ecollement du lieu de branchement}\label{decollement1}
Le but de ce paragraphe est de construire des plongements lisses des sol\'eno\"{\i}des de Sullivan $\mathcal S_g$ 
dans une vari\'et\'e orient\'ee de dimension $4$,
\`a partir de plongements lisses des surfaces branch\'ees $\overline{\Sigma_g}$ (Th\'eor\`eme \ref{decollement}). 

\begin{figure}[h!]\label{pc}
\begin{center}
\input{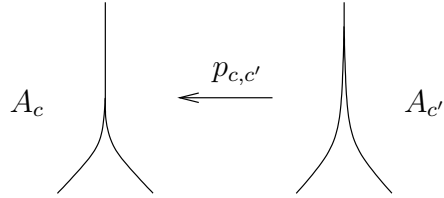}
\caption{L'application $p_{c,c'}$}
\end{center}
\end{figure}

Esquissons la d\'emonstration dans le cas de genre nul, qui est simple et instructive. 
Pour tout $0<c<1$, notons $\overline{A_c}$ le quotient de l'anneau \[A_c=\{ c^2\leq |z|\leq c\}\] par la relation 
qui identifie les points $z$, $-z$ et $z^2$ si $|z|=c$. 
La surface branch\'ee $\overline{A_c}$, hom\'eomorphe \`a $\overline{\Sigma_0}$, est munie de la structure lisse construite en \ref{structure lisse}. 
Les applications  
\[   p_c:O\in \mathcal S \mapsto O\cap A_c \in \overline{A_c},\]
sont les analogues continues des applications $p_{0,n}$. 

Lorsque $1>c'>c>0$, il existe
des applications lisses \[ p_{c,c'}:\overline{A_{c'}}\rightarrow \overline{A_c},\]  
telles que $p_c= p_{c,c'} \circ p_c$. Si $z$ est un point de $A_{c'}$, l'orbite positive 
$O_+(z)= \{ z,r(z),\ldots,r^n(z),\ldots\}$ de $z$ intersecte
l'anneau $A_{c}$ en un point $w$ ou en deux points $w$ et $w^2$ de modules $c$ et $c^2$~: $p_{c,c'}(z)$ 
est la classe de $w$ dans $\overline{A_c}$. 
Lorsque $c'$ est proche de $c$, les fibres de $p_{c,c'}$ sont form\'ees d'un point  
sur l'anneau $\{ c^2 \leq |z|\leq c'^2 \}$, et de deux points sur l'anneau $\{ c'^2\leq |z| \leq c\}$. 
La Figure \ref{pc} repr\'esente l'application $p_{c,c'}$
lorsque $c'$ est proche de $c$. 

Supposons que l'on ait un plongement lisse $\pi_c : \overline{A_c}\rightarrow V$. Choisissons un r\'eel $c'>c$ 
et d\'ecollons dans $V$ les feuillets $D(\sqrt{z})$ 
et $D(-\sqrt{z})$ de $\overline{A_c}$ jusqu'\`a $\{ |z|= c'^2\}$. Nous obtenons alors un plongement 
lisse $\pi_{c'}: \overline{A_{c'}}\rightarrow V$.
D\'ecollons les feuillets d'une distance de plus en plus petite lorsque $c'$ tend vers $1$~; 
la limite de Hausdorff des images des applications $\pi_{c'}$ est un sol\'eno\"{\i}de par surfaces lisses 
dans $V$, qui est diff\'eomorphe au sol\'eno\"{\i}de de Sullivan.

\begin{figure}[h!]\label{decollement.fig}
\begin{center} 
\input{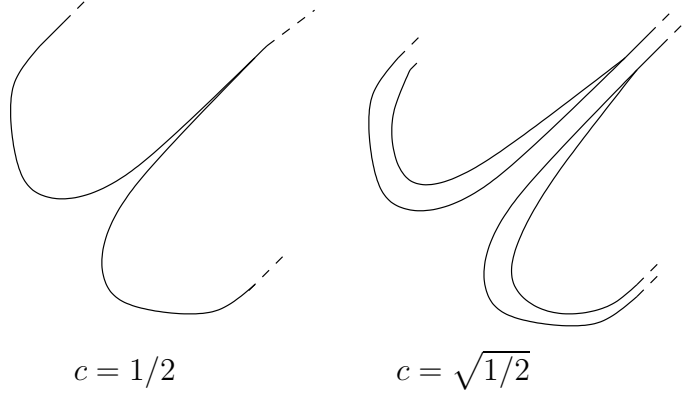}
\caption{D\'ecollement des feuillets dans $M$}
\end{center}
\end{figure}

\begin{theoreme}\label{decollement} \'Etant donn\'ee une vari\'et\'e orient\'ee $V^4$ de dimension $4$ et un plongement 
lisse $\pi : \overline{\Sigma_g}\rightarrow V$, il existe un plongement lisse $\pi_{\infty}:{\mathcal S}_g\rightarrow M$, 
aussi proche que l'on veut de $\pi \circ p_{g,0}$ dans la topologie lisse. \end{theoreme}

\textit{D\'emonstration.} Pour tout entier positif $g$, nous avons un rev\^etement double lisse 
\[ p_{2g,g} : \overline{\Sigma_{2g}} \rightarrow \overline{\Sigma_g},\]
induit par le rev\^etement double $\Sigma_{2g} \rightarrow \Sigma_g$. Le Lemme \ref{plongement d'un revetement} nous explique rigoureusement
le proc\'ed\'e de d\'ecollement des feuillets de $\overline{A_c}$ jusqu'\`a $c'=\sqrt{c}$, et une construction analogue 
en genre sup\'erieur. 

\begin{lemme}\label{plongement d'un revetement} Soit $\pi _g: \overline{\Sigma_g} \rightarrow V$ un plongement lisse  
dans une vari\'et\'e orient\'ee de dimension $4$. Pour tout entier positif $k$,
il existe un plongement lisse $\pi _{2g} : \overline{\Sigma_{2g}} \rightarrow V$ qui est aussi proche que l'on veut 
de l'application lisse $\pi_g\circ p_g$ dans la topologie $C^k$.
(L'application $\pi_{2g}$ est le ``d\'ecollement" du lieu de branchement sur $\Sigma_g$.)  \end{lemme}

\textit{D\'emonstration.} Nous utilisons le Th\'eor\`eme \ref{voisinage global}. 
Il suffit donc de le faire pour les exemples de surfaces 
branch\'ees $\overline{\Sigma}$ plong\'ees dans $V_p$ pour tout $p$. 
Remarquons que $p_{2g,g}:\overline{\Sigma_{2g}}\rightarrow \overline{\Sigma_g}\simeq \overline{\Sigma}$ 
est le rev\^etement double de $\overline{\Sigma}$ d\'efini par les fonctions multivalu\'ees
\[ \sqrt{z}:(A_-\times {\bf C}) \cap \overline{\Sigma}\rightarrow {\bf C}^*,
\ \ \sqrt{z}: (\Sigma_g\times {\bf D}^2)\cap \overline{\Sigma}\rightarrow {\bf C}^*,\]
\[ z^{1/4}: (A_+\times {\bf C})\cap \overline{\Sigma}\rightarrow {\bf C}^*,\]
c'est \`a dire que ces fonctions sont bien d\'efinies sur $\overline{\Sigma_{2g}}$. 
Nous construisons le d\'ecollement en trois parties. Soit $\varepsilon >0$. 

\begin{figure}[h!]\label{pi2g}
\begin{center}
\input{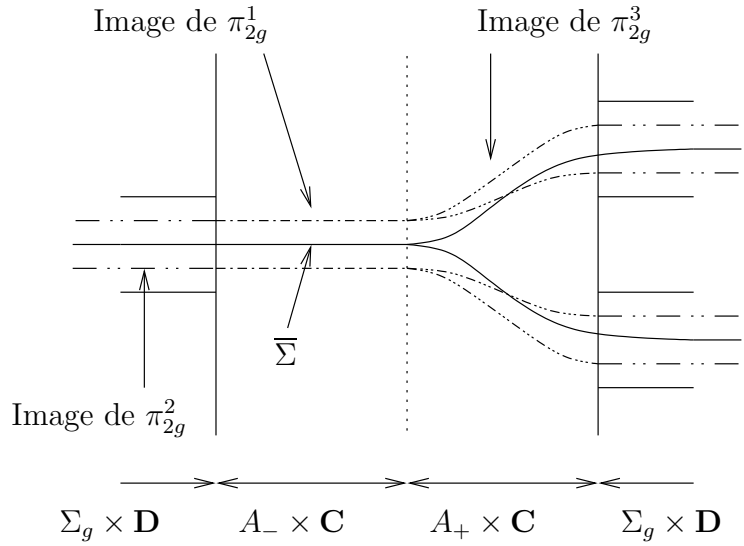}
\caption{ L'application $\pi_{2g}:\overline{\Sigma_{2g}}\rightarrow V_p$ }
\end{center}
\end{figure}

{\bf I.}  \textit{Sur $(A_-\times {\bf C})\cap \overline{\Sigma}$.} Consid\'erons l'application 
\[ \pi_{2g}^{1}: x\in p_{2g,g}^{-1}((A_-\times {\bf C}) 
\cap \overline{\Sigma})\mapsto (x, \varepsilon \sqrt{z})\in A_- \times {\bf C}.\]
Son image est le d\'ecollement du lieu de branchement sur $(A_-\times {\bf C})\cap \overline{\Sigma}$.

{\bf II.} \textit{Sur $\Sigma$.} Soit $\varphi_{\varepsilon} :\Sigma_g \rightarrow (0,+\infty)$ une fonction lisse 
ne s'annulant pas, valant $\varepsilon $ dans un voisinage de $\partial_-\Sigma_g$
et une constante $0<\eta<\varepsilon$ dans un voisinage de $\partial_+\Sigma_{g}$ 
(cette constante sera d\'etermin\'ee plus tard), et telle 
que 
\[  |\varphi_{\varepsilon} |_{C^k} \leq C\varepsilon,\]
la constante $C$ ne d\'ependant pas de $\varepsilon$. Posons  
\[ \pi_{2g}^2: x\in p_{2g,g}^{-1}((\Sigma_g\times {\bf D}^2)\cap \overline{\Sigma}) \mapsto 
(x, \varphi_{\varepsilon}(x) \sqrt{z})\in \Sigma_g\times {\bf D}.\]
L'image de $\pi_{2g}^2$ est le d\'ecollement du lieu de branchement sur $\Sigma$. 

{\bf III.} \textit{Sur $(A_+\times {\bf C})\times \overline{\Sigma}$.} Soit $f_{\varepsilon}:A\rightarrow {\bf C}$ une fonction lisse,  
ne s'annulant pas,
valant $\varepsilon ^2z$ sur $A_-$, et v\'erifiant 
\[  |f^{\frac{1}{2}}-f_{\varepsilon}^{\frac{1}{2}}|_{C^k}\leq \delta_k(\varepsilon),\]
o\`u $\delta_k$ est tend vers $0$ lorsque $\varepsilon $ tend vers $0$. 
D'autre part, soit $g:A\rightarrow {\bf C}$ une fonction lisse, nulle sur $A_-$, ne s'annulant pas sur $A_+$,
valant $z$ sur un voisinage de $\{1\}\times {\bf S}^1$ et telle que les d\'eterminations de $g^{1/4}$ au voisinage de 
$\{ 1/2\}\times {\bf S}^1$ soient aussi lisses. Posons 
\[ \pi_{2g}^3: x\in p_{2g,g}^{-1}((A_+\times {\bf C})\cap \overline{\Sigma}) 
\mapsto (x, f_{\varepsilon}^{\frac{1}{2}} +\eta g^{\frac{2p+1}{4}}).\]
Nous choisissons $\eta$ assez petit en sorte que $\pi_{2g}^3$ soit un plongement lisse. 
L'image de $\pi_{2g}^3$ est le d\'ecollement du lieu de 
branchement sur $(A_+\times {\bf C} )\cap \overline{\Sigma}$.

Les applications $\pi_{2g}^i,\ i=1,2,3$ 
se recollent bien avec les changements de cartes $\Phi_+$ et $\Phi_-$, pour former un plongement  
lisse $\pi_{2g}:\overline{\Sigma_{2g}}\rightarrow V_p $ (voir Figure), v\'erifiant 
\[d_k(\pi_{2g},\pi_g\circ p_{2g,g}) \leq \sup(\delta(\varepsilon),C\varepsilon),\]
et le Lemme est d\'emontr\'e. 

\begin{remarque}\label{remarque} Au lieu de ``d\'ecoller" le lieu de branchement jusqu'\`a $\partial_+\Sigma$ dans la deuxi\`eme \'etape, 
on peut s'arr\^eter sur une courbe ferm\'ee simple contenue dans l'int\'erieur de $\Sigma$, s\'eparant 
$\partial_-\Sigma$ de $\partial_+\Sigma$. Nous obtenons de cette fa\c{c}on une surface branch\'ee lisse $\overline{\Sigma}'$ 
de genre $g\leq g(\overline{\Sigma}')\leq 2g$ donn\'e, arbitrairement proche de $\overline{\Sigma}$.\end{remarque}

Achevons la d\'emonstration du Th\'eor\`eme \ref{decollement}. Soit $\overline{V}$ un voisinage compact de la surface 
branch\'ee $\pi(\overline{\Sigma_g})$. Pour tout $k$, notons $C^{k}(\mathcal S_g,\overline{V})$ l'espace des applications de
classes $C^k$ de $\mathcal S_g$ dans $\overline{V}$, et $d_k$ la distance sur $C^k({\mathcal S_g},\overline{V})$ d\'efinie par 
\[  d_k (f,g)= \sup_{x\in \mathcal S_g} d(j^k(f)(x),\ j^k(g)(x)),\]
o\`u $j^k$ d\'esigne le $k$-i\`eme jet, et la distance $d$ sur les jets \'etant induite par des m\'etriques donn\'ees sur $\mathcal S_g$
et $\overline{V}$. Pour tout $k$, les espaces $(C^k({\mathcal S_g},\overline{V}),d_k)$ sont complets. 

Pour tout $n$, notons $\varepsilon _n>0$ le supremum du diam\`etre des fibres de $p_{g,n}$. La suite $\varepsilon _n$ 
tend vers $0$, d'apr\`es le Lemme \ref{approximation}. 
En appliquant le Lemme pr\'ec\'edent, nous construisons une suite de plongements \textit{lisses} 
\[ \pi_{n} :\overline{\Sigma_{2^ng}}\rightarrow {\bf C}P^2,\] 
et des voisinages $\mathcal V_n$ de $\pi_{n}\circ p_{g,n}\in C^k(\mathcal S_g)$ pour la topologie $C^k$ tels que~:

\ \ (i)  Tout $\pi'\in \mathcal V_n$ est une immersion $2\varepsilon _n$-injective, i.e. deux points s\'epar\'es d'une distance 
sup\'erieure \`a $2\varepsilon _n$ n'ont pas la m\^eme image par $\pi'$. 

\ (ii) Pour tout $n'\geq n$, $\pi_{n'} \in \mathcal V_n$. 

(iii) Pour tout $k$, le $d_k$-diam\`etre de $\mathcal V_n$ tend vers $0$ lorsque $n$ tend vers l'infini.\\
\noindent La suite $\pi_n \circ p_{g,n}$ tend alors vers un plongement lisse de $\mathcal V_0$ et le Th\'eor\`eme 
\ref{decollement} est d\'emontr\'e.

\begin{question} Le lecteur pourra s'assurer que la tresse normale de la 
surface branch\'ee $\overline{\Sigma'}= \pi_{2g}(\overline{\Sigma_{2g}})$ construite au Lemme \ref{plongement d'un revetement} 
est la m\^eme que celle de la surface $\overline{\Sigma}$. Est-il possible de d\'efinir la tresse normale de
l'image du sol\'eno\"{\i}de $\mathcal S_g$ par un plongement lisse \`a valeurs dans une vari\'et\'e orient\'ee 
de dimension $4$?  \end{question}

\section{Surfaces branch\'ees et sol\'eno\"ides symplectiques plans}\label{surfaces branchees symplectiques}

Dans ce paragraphe, nous construisons des plongements lisses $\varepsilon$-holomorphes de 
la surface branch\'ee $\overline{\Sigma}_2$ dans le plan projectif complexe, 
et donc du sol\'eno\"{\i}de $\mathcal S _2$. 
Le r\'esultat d\'ecoule de consid\'erations faciles sur la topologie de certaines 
courbes alg\'ebriques planes. 

\vspace{0.2cm}

Pour tout $r>0$, soit $B_r =\{ |x| \leq r,\ |y|\leq r\} \subset {\bf C}^2$. 

\vspace{0.2cm}

{\bf Fait I.} \textit{Dans $B_2$, consid\'erons la famille de courbes holomorphes $C_{\varepsilon}$ d\'efinies 
par les \'equations 
\[  xy^2 = \varepsilon ^2 + \varepsilon ^3R(x),\]
o\`u $R(x)$ est une fonction holomorphe d\'efinie sur un voisinage du disque ferm\'e de rayon $2$. 
Pour tout $\varepsilon>0$ suffisamment petit, il existe un plongement lisse $c\varepsilon$-holomorphe de $\overline{A}$ dans  
$B_2$, dont l'image $\overline{A}_{\varepsilon}$ co\"{\i}ncide avec $C_{\varepsilon} $ sur $B_2-B_{2-1/3}$, 
o\`u $c>0$ est une constante ne d\'ependant que de $R$.\label{chirurgie}}

\vspace{0.2cm} 
 
\textit{D\'emonstration.} Notons $B_x=\{1\leq |x| \leq 2, \ |y|\leq 2\}$ et $B_y =\{ |x|\leq 2,\ 1\leq |y| \leq 2\}$. 
Nous avons $B_2 = B_x \cup B_1 \cup B_y$.

Consid\'erons une fonction lisse $\rho_1 : [1,2] \rightarrow {\bf R}$ 
qui vaut $1$ sur $[2-1/3,2]$, $0$ sur $[1,1+1/3]$ et telle que $0<\rho_1<1$ sur $[1+1/3,2-1/3]$. D\'efinissons 
\[ A_1 := \{ (x,y)\in B_x \  |\ xy^2 = \rho_1(|x|)  (\varepsilon ^2 + \varepsilon ^3 R(x)) \}.\]
Observons que lorsque $\varepsilon>0$ est assez petit, $\overline{A}_1$ n'intersecte pas $B_x \cap B_y$. De plus,
$\overline{A}_1$ est un plongement lisse de $\overline{A}$.

Consid\'erons une fonction lisse $\rho_2 :B_1\rightarrow {\bf R}$ 
qui vaut $1$ sur $B_{1/3}$, $0$ \`a l'ext\'erieur de $B_{2/3}$ et qui prend des valeurs strictement positives sur $Int(B_{2/3})-B_{1/3}$. 
D\'efinissons 
\[  A_2 := \{ (x,y)\in B_1 \  |\ xy = \rho (x,y) \varepsilon\}\ ;\]
$A_2$ est un anneau lisse, et l'union $A_1 \cup A_2$ est l'image d'un plongement lisse de $\overline {A}$.

Par le Th\'eor\`eme des fonctions implicites, il existe une
famille de fonctions holomorphes $X_{\varepsilon}$ d\'efinies sur l'anneau $\{ 1\leq |y| \leq 2\}$ telles que 
\begin{itemize} 
\item $X_{\varepsilon} = O_{\varepsilon\rightarrow 0}(\varepsilon)$.  
\item $C_{\varepsilon} \cap B_y = \{ (x,y)\ | \ x = X_{\varepsilon} (y)\}. $ 
\end{itemize} 
D\'efinissons 
\[  A_3 :=  \{ (x,y)\in B_y \  |\ x= \rho_1(|y|) X_{\varepsilon} (y)\}.\]
\`A nouveau, $A_3$ est un anneau, et l'union $A_1 \cup A_2 \cup A_3$ est l'image par un plongement lisse de $\overline{A}$,
qui co\"incide avec $C_{\varepsilon}$ sur $B_2 - B_{2-1/3}$. 

Il est ais\'e de v\'erifier que $A_1 \cup A_2 \cup A_3$ est une surface branch\'ee $c\varepsilon$-holomorphe, 
o\`u $c>0$ est une constante ne d\'ependant que de $R$. 
Le Fait I est d\'emontr\'e. 

\vspace{0.2cm} 

Consid\'erons la famille de 
courbes alg\'ebriques $C_{\varepsilon}$ de ${\bf C}P^2$ 
d\'efinies par les \'equations homog\`enes 
\[  xy^2z - \varepsilon ^2 z^4 - \varepsilon ^3 R(x,z) =0,\]
o\`u $R$ est un polyn\^ome homog\`ene de degr\'e $4$ en $x$ et $z$. Si $R$ est choisi g\'en\'eriquement,
alors pour tout $\varepsilon > 0$ la courbe $C_{\varepsilon}$ est non singuli\`ere. 
Nous fixerons l'un de ces choix de $R$. 
\'Etant de degr\'e $4$ et non singuli\`eres, les courbes $C_{\varepsilon}$ sont de genre $3$ lorsque $\varepsilon >0$. 

Dans une surface, un anneau est dit \textit{essentiel} s'il ne s\'epare pas la surface. 

\vspace{0.2cm}
 
{\bf Fait II.} \textit{Pour $\varepsilon>0$ assez petit, $B_2\cap C_{\varepsilon}$ est un anneau essentiel dans $C_{\varepsilon}$.}

\vspace{0.2cm}

\textit{D\'emonstration.} Dans le plan affine ${\bf C}^2$ de coordonn\'ees $x$ et $y$, l'\'equa\-tion de $C_{\varepsilon}$
est 
\[  xy^2 = \varepsilon^2 + \varepsilon  ^3 R(x,1).\] 
Soient 
\[ D=\{ x\in {\bf C} \ |\ |x|\leq 2\} \ \ \mathrm{et}\ \ 
A=\{ |x|\leq 2,\ |\varepsilon ^2 +\varepsilon ^3 R(x,1)|\leq 4|x|\}\subset D.\] 
Lorsque $\varepsilon$ est assez petit $A$ est un anneau. 
Comme la projection $x:B_2\cap C_{\varepsilon} \rightarrow A$ 
est un rev\^etement double de $C_{\varepsilon} \cap B_2$ sur $A$,  
$C_{\varepsilon} \cap B_2$ est un anneau. 

Pour voir que $C_{\varepsilon}-B_2$ est connexe pour $\varepsilon >0$ petit, 
il suffit de le v\'erifier pour $\varepsilon =0$, ce qui est trivial. Le Fait II est d\'emontr\'e. 

\begin{theoreme} Pour tout $\varepsilon>0$, il existe un plongement lisse $\varepsilon$-holomorphe de $\overline{\Sigma}_2$ dans 
le plan projectif complexe.\label{plgt surface branchee}\end{theoreme} 

\textit{D\'emonstration.} Il suffit de regarder la famille de courbes $C_{\varepsilon}$ 
\'etudi\'ee au Fait II, et de lui appliquer la chirurgie expliqu\'ee par le Fait I dans la boule $B_2$. 
On obtient une surface branch\'ee $c\varepsilon$-holomorphe $\overline{\Sigma}_{\varepsilon}$, 
obtenue en enlevant \`a une surface de genre $3$ un anneau essentiel et en le rempla\c{c}ant par l'anneau branch\'e $\overline{A}$~:  
$\overline{\Sigma}_{\varepsilon}$ est donc l'image d'un plongement lisse $c\varepsilon$-holomorphe 
de $\overline{\Sigma}_2$ et le Th\'eor\`eme est d\'emontr\'e.

\begin{corollaire}  Pour tout $\varepsilon>0$, il existe un plongement lisse $\varepsilon$-holomorphe de $ {\mathcal S}_2$
dans le plan projectif complexe.\label{plgt solenoide}\end{corollaire}

\textit{D\'emonstration.} Ceci d\'ecoule des Th\'eor\`emes \ref{plgt surface branchee} et \ref{decollement}.

\vspace{0.2cm}

Nous terminons par deux remarques.

\begin{remarque}\label{conjecture} Cet exemple de sol\'eno\"{\i}de $\varepsilon$-holomorphe dans ${\bf C}P^2$ 
sem\-ble indiquer qu'il existe un sol\'eno\"{\i}de de cette nature dont les feuilles sont des courbes holomorphes. 
Il est fond\'e sur le principe, d\^u \`a Gromov, selon lequel une surface symplectique est isotope \`a une courbe 
holomorphe~\cite{Gromov1}. 
L'approche consiste \`a utiliser la th\'eorie des courbes $J$-holomor\-phes pour des structures 
presque-complexes $J$ calibr\'ees sur la forme de Fubini-Study. 
Dans notre cas, il est possible de construire une structure presque-complexe $J_0$ calibr\'ee sur la forme de Fubini-Study 
pour laquelle les feuilles du sol\'eno\"{\i}de de \ref{plgt solenoide} soient $J_0$-holomorphes. 
Comme l'espace des structures presque-complexes calibr\'ees est connexe par arcs,
nous pouvons choisir une famille continue $J_t$ de structures complexes liant $J_0$ \`a la structure standard $J_1$. Il s'agit alors de 
d\'eformer notre sol\'eno\"{\i}de initial en un sol\'eno\"{\i}de par courbes $J_t$-holomorphes pour tout $t$. 
\`A titre de motivation nous sugg\'erons les r\'ef\'erences suivantes o\`u des r\'esultats 
semblables ont \'et\'e d\'emontr\'es~:~\cite{Gromov,Moser}. M\^eme si 
ces techniques de d\'eformation semblent hors de port\'ee pour l'instant, nous sommes amen\'es \`a formuler la conjecture suivante.  

\begin{conjecture} Il existe un sol\'eno\"{\i}de compact par courbes holomorphes dans le plan projectif complexe, ne 
contenant pas de courbe alg\'ebrique. Nous n'avons par contre pas d'intuition sur la question de savoir si ce sol\'eno\"{\i}de 
est tangent \`a un feuilletage holomorphe. \end{conjecture}
\end{remarque}

\begin{remarque} [A. Glutsyuk] Il est tentant de penser que puisque ${\mathcal S}_2 \subset {\bf C}P^2$ n'a pas 
d'homologie en dimension $2$, il peut \^etre isotop\'e symplectiquement dans ${\bf R}^4$, muni de la forme symplectique standard. 
Pourtant, comme nous l'a communiqu\'e A. Glutsyuk ce n'est pas le cas. En effet, il existerait alors une structure presque complexe 
sur ${\bf C}P^2$ compatible avec la structure symplectique, pour laquelle il existe un sol\'eno\"{i}de compact $J$-holomorphe 
n'intersectant pas la droite $J$-holomorphe de l'infini. Or d'apr\`es les travaux de Gromov \cite{Gromov1}, 
par deux 
points distincts passent une unique courbe rationnelle $J$-holomorphe homologue \`a la droite de l'infini.
Il y aurait donc au moins deux composantes connexes dans l'espace des courbes $J$-holomorphes homologue \`a la droite de l'infini~:
celles qui coupent le sol\'eno\"{\i}de et celles qui ne le coupent pas (car les intersections sont positives).  
Ceci est une contradiction.  \end{remarque}

\end{document}